\documentclass{amsart}

\usepackage{fullpage}
\usepackage{amsmath}
\usepackage{amsthm}
\usepackage{amsfonts}
\usepackage{amssymb}

\input xy
\xyoption{matrix} \xyoption{arrow} \xyoption{curve}

\begin{document}


\newcommand{\Hom}{\mathrm{Hom}}
\newcommand{\RHom}{\mathrm{RHom}^*}
\newcommand{\HOM}{\mathrm{HOM}}
\newcommand{\stHom}{\underline{\mathrm{Hom}}}
\newcommand{\Ext}{\mathrm{Ext}}
\newcommand{\Tor}{\mathrm{Tor}}
\newcommand{\HH}{\mathrm{HH}}
\newcommand{\Endo}{\mathrm{End}}
\newcommand{\ENDO}{\mathrm{END}}
\newcommand{\stEndo}{\mathrm{\underline{End}}}
\newcommand{\Tr}{\mathrm{Tr}}
\newcommand{\spec}{\mathrm{Spec}}
\newcommand{\GL}{\mathrm{GL}}


\newcommand{\coker}{\mathrm{coker}}
\newcommand{\aut}{\mathrm{Aut}}
\newcommand{\out}{\mathrm{Out}}
\newcommand{\op}{\mathrm{op}}
\newcommand{\add}{\mathrm{add}}
\newcommand{\ADD}{\mathrm{ADD}}
\newcommand{\ind}{\mathrm{ind}}
\newcommand{\rad}{\mathrm{rad}}
\newcommand{\soc}{\mathrm{soc}}
\newcommand{\ann}{\mathrm{ann}}
\newcommand{\im}{\mathrm{im}}
\newcommand{\chr}{\mathrm{char}}
\newcommand{\pdim}{\mathrm{p.dim}}
\newcommand{\dgn}{\leq_{\mathrm{deg}}}
\newcommand{\cx}{\mbox{cx}}


\newcommand{\rmod}{\mbox{mod-}}
\newcommand{\Rmod}{\mbox{Mod-}}
\newcommand{\lmod}{\mbox{-mod}}
\newcommand{\lMod}{\mbox{-Mod}}
\newcommand{\stmod}{\mbox{\underline{mod}-}}
\newcommand{\stlmod}{\mbox{-\underline{mod}}}
\newcommand{\Mod}{\mathcal{M}\mbox{od}}

\newcommand{\gmod}[1]{\mbox{mod}_{#1}\mbox{-}}
\newcommand{\gMod}[1]{\mbox{Mod}_{#1}\mbox{-}}
\newcommand{\Bimod}[1]{\mathrm{Bimod}_{#1}\mbox{-}}

\newcommand{\proj}{\mbox{proj-}}
\newcommand{\lproj}{\mbox{-proj}}
\newcommand{\Proj}{\mbox{Proj-}}
\newcommand{\inj}{\mbox{inj-}}
\newcommand{\coh}{\mbox{coh-}}


\newcommand{\und}[1]{\underline{#1}}
\newcommand{\gen}[1]{\langle #1 \rangle}
\newcommand{\floor}[1]{\lfloor #1 \rfloor}
\newcommand{\ceil}[1]{\lceil #1 \rceil}
\newcommand{\bnc}[2]{\left( \begin{array}{c} #1 \\ #2 \end{array} \right) }
\newcommand{\bimo}[1]{{}_{#1}#1_{#1}}
\newcommand{\ses}[5]{\ensuremath{0 \rightarrow #1 \stackrel{#4}{\longrightarrow} 
#2 \stackrel{#5}{\longrightarrow} #3 \rightarrow 0}}
\newcommand{\A}{\mathcal{A}}
\newcommand{\B}{\mathcal{B}}
\newcommand{\uB}{\underline{\mathcal{B}}}
\newcommand{\C}{\mathcal{C}}
\newcommand{\Or}{\mathcal{O}}
\newcommand{\uC}{\underline{\mathcal{C}}}
\newcommand{\D}{\mathcal{D}}
\newcommand{\E}{\mathcal{E}}
\newcommand{\F}{\mathcal{F}}
\newcommand{\p}{\mathcal{P}}
\newcommand{\m}{\mathfrak{m}}
\newcommand{\ul}[1]{\underline{#1}}


\newtheorem{theorem}{Theorem}
\newtheorem{defin}[theorem]{Definition}
\newtheorem{propos}[theorem]{Proposition}
\newtheorem{lemma}[theorem]{Lemma}
\newtheorem{coro}[theorem]{Corollary}
\newtheorem{qu}[theorem]{Question}

\title{On periodicity in bounded projective resolutions}
\date{\today}
\author{Alex Dugas}
\address{Department of Mathematics, University of the Pacific, 3601 Pacific Ave, Stockton, CA 95211, USA}
\email{adugas@pacific.edu}

\subjclass[2010]{16E05}
\keywords{periodic resolution, complexity, open orbit, rigid module, stable equivalence of Morita type}

\begin{abstract} Let $A$ be a self-injective algebra over an algebraically closed field $k$.  We show that if an $A$-module $M$ of complexity one has an open $\GL_d(k)$-orbit in the variety $\Mod^A_d$ of $d$-dimensional $A$-modules, then $M$ is periodic.  As a corollary we see that any simple $A$-module of complexity one must be periodic.   In the course of the proof, we also show that modules with open orbits are preserved by stable equivalences of Morita type between self-injective algebras.
\end{abstract}

\maketitle
 
\setcounter{equation}{0}

The complexity of a module $M$ over a finite-dimensional algebra $A$ is an invariant measuring the growth rate of a minimal projective resolution $P^{\bullet}$ of $M$:
$$\cx_A(M) := \inf\ \{n \geq 0\ |\ \exists\ C>0\ \mbox{s.t.}\ \dim_k\ P^t \leq C t^{n-1}\ \forall\ t \geq 0 \}.$$ 
For example, $\cx_A(M) = 0$ if and only if $M$ has finite projective dimension, and $\cx_A(M) \leq 1$ if and only if $M$ has a projective resolution with terms of bounded dimension.  The study of complexity emerged in the representation theory of finite groups, motivated in part by results of Alperin \cite{Alp} and Eisenbud  \cite{Eis}, stating that any module of complexity one over a group algebra $kG$ is {\it periodic}, i.e., has a periodic projective resolution.  However, the same is not true over arbitrary self-injective algebras; the simplest counterexamples being afforded by the quantum exterior algebras in two generators \cite{LiuSchulz}.  Namely, if $\Lambda_q = k\gen{x,y}/(x^2, y^2, xy+qyx)$ for $q \in k^*$ of infinite multiplicative order, and $M = \Lambda_q (x+y)$, then it is not hard to see that $\Omega^iM \cong \Lambda_q (x+q^i y)$ for all $i \in \mathbb{Z}$.  Moreover, these syzygies are all $2$-dimensional and pairwise non-isomorphic.  Similar counterexamples over commutative local rings were first given by Gasharov and Peeva \cite{GashPeev}, and now many such are known.  Our goal here is to show that any simple module of complexity one over a self-injective algebra is indeed periodic.  

We establish this result through a geometric argument using the affine varieties $\Mod^A_d$ parametrizing $d$-dimensional $A$-modules.  The general linear group $\GL_d(k)$ acts morphically on $\Mod^A_d$ via conjugation, and the $\GL_d(k)$-orbits are in one-to-one correspondence with the isomorphism classes of $d$-dimensional $A$-modules.  If $M$ is a $d$-dimensional $A$-module, we identify $M$ with the corresponding point in $\Mod^A_d$, and we write $\Or(M)$ for the $\GL_d(k)$-orbit of $M$.  We will give an algebraic characterization for $\Or(M)$ to be Zariski open, which we can use to deduce the following key result.
  
\begin{theorem}  Let $A$ be a basic self-injective algebra over an algebraically closed field $k$ and $M_A$ a $d$-dimensional $A$-module.  If $M$ has complexity one and an open $\GL_d(k)$-orbit in $\Mod^A_d$, then $M$ is periodic.
\end{theorem}

When $A$ is basic, any simple $A$-module is one-dimensional, and hence its orbit in $\Mod^A_1$ is trivially open.  Since periodicity of a module is clearly preserved by Morita equivalence, we obtain our previous claim:

\begin{coro} Any simple module of complexity one over a self-injective algebra is periodic.
\end{coro}

In fact, these two results are essentially un-noted corollaries of work of Everett Dade \cite{Dade}, which we discovered independently and by a less direct approach making use of more recent developments.  While our definitions differ somewhat from those of Dade, our preliminary results closely parellel his own, as noted below.

\section{Open orbits and rigidity}
We assume that $A$ is a finitely-generated algebra over an algebraically closed field $k$.  Fix a presentation of the algebra $A$ by generators and relations $A \cong k\gen{x_1, \ldots, x_s}/(\rho_i)_{i \in I}$ with $\rho_{i_0} = x_1-1$.  Clearly, a $d$-dimensional unital (right) $A$-module corresponds to a choice of $s$ $d \times d$ matrices $X_1 = I , X_2, \ldots, X_s$ satisfying the relations $\rho_i$ for each $i \in I$.  As each such relation can be expressed via polynomial equations in the entries of the matrices $X_i$, we see that the possible $A$-module structures on $k^d$ are parametrized by the points of the (usually reducible) affine variety
$$\Mod^A_d = \{(X_1, \ldots, X_s) \in M_n(k)^d\ |\ X_1= I,\ \rho_i(X_1,\ldots,X_s)=0\ \forall\ i \in I\}.$$
The isomorphism class of the module corresponding to a point in this variety obviously does not depend on a fixed basis of $k^d$ and is thus independent of a change of bases.  Indeed, the general linear group $\GL_d(k)$ acts morphically on this variety via conjugation (corresponding to rewriting the matrices $X_i$ with respect to a different basis of $k^d$), and the orbits of this action correspond precisely to the isomorphism classes of $d$-dimensional $A$-modules.   If $M$ is a $d$-dimensional $A$-module, with a slight abuse of notation, we identify $M$ with the corresponding point in $\Mod^A_d$ and write $\Or(M)$ for the $\GL_d(k)$-orbit of $M$.    Since any orbit $\Or(M)$ is open in its closure (i.e., locally closed) and the irreducible components of $\Mod^A_d$ are $\GL_d(k)$-stable, it is not hard to see that $\Or(M)$ is open if and only if its closure $\overline{\Or(M)}$ is an irreducible component of $\Mod^A_d$.  Finally, we recall that for two $d$-dimensional $A$-modules $M$ and $N$, $M$ is said to {\it degenerate} to $N$ if $N \in \overline{\Or(M)}$, and we express this using the notation $M \dgn N$.

Our present goal is to formulate an algebraic equivalent of openness of the orbit $\Or(M)$ in $\Mod^A_d$.  Indeed, we will see that it is captured by the following definition.  As usual, we write $\spec\ R$ for the maximal ideal spectrum of  a commutative noetherian $k$-algebra $R$.

\begin{defin}  We say that a $d$-dimensional $A$-module $M_A$ is {\bf rigid} if for every one-dimensional noetherian domain $R$ (over $k$) and all $(R,A)$-bimodules ${}_R \tilde{M}_A$ which are free of rank $d$ over $R$, the set $$\{\m \in \spec\ R\ |\ R/\m \otimes_R \tilde{M}_A \cong M_A\}$$ is open in $\spec\ R$.  (Note that since $\spec\ R$ is one-dimensional, it would be equivalent to require this subset to be dense whenever it is non-empty.)
\end{defin}

Geometrically, this definition says $M$ is rigid if and only if  $\varphi^{-1}(\Or(M))$ is open for all morphisms $\varphi$ from an affine curve $\C$ to $\Mod^A_d$.  To see this, we briefly consider the affine scheme $\Mod^A_d(-)$,  viewing it as a functor from $k$-algebras to sets, which is represented by the coordinate ring $\Gamma$ of $\Mod^A_d$ as in \cite{Zwa}.  An $R$-valued point of this scheme is hence given by a $k$-algebra homomorphism $\Gamma \rightarrow R$, which corresponds to placing a (right) $A$-module structure on $R^d$.  From this point of view, a point $u$ of $\Mod^A_d(R)$ corresponds to an $(R,A)$-bimodule $U$, free of rank $d$ over $R$.  Moreover, if $\phi : R \rightarrow S$ is a $k$-algebra homomorphism, the $(S,A)$-bimodule corresponding to the point $\Mod^A_d(\phi)(u) \in \Mod^A_d(S)$ is given by $S \otimes_R U$, where the $R$-action on $S$ is induced by $\phi$.

  Now let $\C = \spec\ R$ be an irreducible affine curve and $\varphi : \C \rightarrow \Mod^A_d$ a morphism with $\varphi^* : \Gamma \rightarrow R$ the induced map of coordinate rings.  We thus obtain an $(R,A)$-bimodule $\tilde{M}$, free of rank $d$ over $R$, corresponding to $\varphi$.  If we identify $\m \in \spec\ R$ with a $k$-algebra homomorphism $f_{\m} : R \rightarrow k$, then we have $\varphi(\m) = f_{\m} \varphi^* = \Mod^A_d(f_{\m})(\varphi^*) \in \Mod^A_d(k)$. Thus the $A$-module corresponding to $\varphi(\m)$ is $k \otimes_R \tilde{M}_A = R/\m \otimes_R \tilde{M}_A$ by the last paragraph.  Consequently we see that the set of $\m \in \spec\ R$ for which $R/\m \otimes_R \tilde{M}_A \cong M_A$ is precisely $\varphi^{-1}(\Or(M))$.

\begin{propos}  A $d$-dimensional $A$-module $M$ is rigid if and only if $\Or(M)$ is open in $\Mod^A_d$.   
\end{propos}

\noindent
{\it Proof.} Let $\C$ be an irreducible affine curve and $\varphi: \C \rightarrow \Mod^A_d$ a morphism.  Assuming that $\Or(M)$ is open in $\Mod^A_d$, we see that $\varphi^{-1}(\Or(M))$ must also be open in $\C$.  Thus rigidity of $M$ follows from the above discussion.

Conversely, suppose $M$ is rigid and consider a module $N$ in the same irreducible component of $\Mod^A_d$ as $M$.  Then $M$ and $N$ can be connected by an irreducible affine curve.  In other words, there is a one-dimensional noetherian domain $R$ and a morphism $\varphi : \spec\ R \rightarrow \Mod^A_d$ with $\varphi(\m_0) = M$ and $\varphi(\m_1) = N$ for some $\m_0, \m_1 \in \spec\ R$.  From the comments above, $\varphi$ corresponds to an $(R,A)$-bimodule ${}_R \tilde{M}_A$, free of rank $d$ over $R$, such that $R/\m_0 \otimes_R \tilde{M} = M$ and $R/\m_1 \otimes_R \tilde{M} = N$.  Since $M$ is rigid, we have $\varphi(\m) \in \Or(M)$ for a dense subset of $\m \in \spec\ R$, and it follows that $N \in \overline{\Or(M)}$ (see II.2.4-2.5 and II.3.5 in \cite{Kraft}).  Thus $\overline{\Or(M)}$ is an irreducible component of $\Mod^A_d$.  $\Box$\\

We end this section with a lemma giving a slightly stronger characterization of rigidity, which we will need later.

\begin{lemma} Let $M_A$ be a rigid $d$-dimensional $A$-module $M_A$.  Then  the set $$\{\m \in \spec\ R\ |\ R/\m \otimes_R \tilde{M}_A \cong M_A\}$$ is open in $\spec\ R$ for every one-dimensional noetherian domain $R$ and all $(R,A)$-bimodules ${}_R \tilde{M}_A$ which are locally free of rank $d$ over $R$.
\end{lemma}

\noindent
{\it Proof.}  Suppose $M_A$ is rigid and let ${}_R \tilde{M}_A$ be an $(R,A)$-bimodule which is locally free of rank $d$ over $R$.  Moreover, suppose that $R/\m_0 \otimes_R \tilde{M}_A \cong M_A$ for some $\m_0 \in \spec\ R$.  According to Exercises 4.11-12 of \cite{CAVG}, we can find $f \in R \setminus \m_0$ such that $\tilde{M}[f^{-1}]$ is free over $R[f^{-1}]$.  Since $R[f^{-1}]$ is again a one-dimensional noetherian domain, and $R[f^{-1}]/\m_0[f^{-1}] \otimes_{R[f^{-1}]} \tilde{M}[f^{-1}] \cong R/\m_0 \otimes_R \tilde{M} \cong M_A$, we have $R[f^{-1}]/\m \otimes_{R[f^{-1}]} \tilde{M}[f^{-1}] = R[f^{-1}]/\m \otimes_R \tilde{M} \cong M_A$ for a dense set of $\m \in \spec\ R[f^{-1}]$.  Since $\spec\ R[f^{-1}]$ can be identified with the open subset of $\spec\ R$ consisting of maximal ideals not containing $f$, we get $R/\m \otimes_R \tilde{M}_A \cong M_A$ for a dense set of $\m \in \spec\ R$.  $\Box$\\

\section{Stable equivalences of Morita type}

We now assume that $A$ and $B$ are self-injective in order to study the effect of a stable equivalence on modules with open orbits.  Of course, if we hope to show that rigidity is preserved by an equivalence of stable categories, then we should expect that rigidity is not affected by adding or removing projective direct summands.

\begin{lemma}[Cf. Theorem 3.16 in \cite{Dade}]  Let $M$ be a $d$-dimensional $A$-module and let $P$ be a projective $A$-module.   Then $\Or(M \oplus P)$ is  open if and only if $\Or(M)$ is open.
\end{lemma}

\noindent
{\it Proof.}  $\Leftarrow$:  This direction follows from Theorem 1.3 in \cite{CBS}.  Suppose $\Or(M)$ is open.  Then the set $\E(\Or(M) \times \Or(P))$, defined as the $\GL$-stable subset consisting of extensions $E$ of $M$ by $P$, coincides with $\Or(M \oplus P)$ and is open by part (iii) of the aforementioned theorem.  
$\Rightarrow$:  Let $C_M$ and $C_P$ be irreducible components containing $M$ and $P$ respectively.  Since $P$ is projective, $\Ext^1_A(N,P) = 0=\Ext^1_A(P,N)$ for all $N \in C_M$, and hence $\overline{C_M \oplus C_P}$ is an irreducible component by Theorem 1.2 in \cite{CBS}.  Since $M \oplus P \in C_M \oplus C_P$ and $\Or(M \oplus P)$ is open by assumption, we must have $\overline{\Or(M \oplus P)} = \overline{C_M \oplus C_P}$.  Thus if $N \in C_M$, we have $M \oplus P \dgn N \oplus P$.  By Bongartz's cancellation result \cite{Bon}, we can cancel the projective $P$ to conclude that $M \dgn N$.  This shows that $\C_M = \overline{\Or(M)}$.  $\Box$ \\

We now recall the definition of a stable equivalence of Morita type.  An $(A,B)$-bimodule $X$ is called {\it left-right projective} if ${}_A X$ and $X_B$ are both projective.

\begin{defin} A pair of left-right projective bimodules ${}_A X_B$ and ${}_B Y_A$ is said to induce a {\bf stable equivalence of Morita type} between $A$ and $B$ if there are bimodule isomorphisms $${}_A X \otimes_B Y_A \cong A \oplus P\ \ \ \mbox{and}\  \ \ {}_B Y \otimes_A X_B \cong B \oplus Q$$ for projective bimodules $P$ and $Q$.  In this case the exact functors $-\otimes_A X_B$ and $-\otimes_B Y_A$ induce inverse equivalences between the stable categories $\stmod A$ and $\stmod B$.
\end{defin}

The main examples of interest to us here are the syzygy functors over a self-injective algebra.  Namely, when $A$ is self-injective, tensoring with the $(A,A)$-bimodules $\Omega^n_{A^e}(A)$ and $\Omega^{-n}_{A^e}(A)$ induces the $n^{th}$ syzygy functor $\Omega^n$ and its inverse $\Omega^{-n}$, respectively, on $\stmod A$.   By Theorem 3.2 of \cite{AMRT}, this pair of bimodules yields a stable equivalence of Morita type.   Likewise, since the auto-equivalences $D\Tr$ and $\Tr D$ of $\stmod A$ are induced by exact functors on $\rmod A$ \cite{TEG}, they also lift to a stable equivalence of Morita type.  Additional examples are found in modular representation theory.  Namely, if $B$ is a block of a finite group $G$ with a trivial-intersection defect group $D$, and $b$ is the Brauer correspondent of $B$, then the Green correspondence (i.e., the induction and restriction functors) provides a stable equivalence of Morita type between $B$ and $b$.  In fact, in this context Dade has shown that the Green correspondence preserves algebraically rigid modules.  Using our definition of rigidity, it is not hard to extend this result to arbitrary stable equivalences of Morita type between self-injective algebras. 

\begin{propos}  [Cf. Theorem 4.4 in \cite{Dade}]  Assume ${}_A X_B$ and ${}_B Y_A$ give a stable equivalence of Morita type between $A$ and $B$.  Then $\Or(M)$ is open in $\Mod^A_d$ if and only if $\Or(M \otimes_A X_B)$ is open in $\Mod^B_e$, where $d = \dim_k M$ and $e = \dim_k (M \otimes_A X_B)$.
\end{propos}

\noindent
{\it Proof.}  $\Leftarrow$: By the lemma, we may assume that $M$ contains no projective summands.  We know that $\Or(M)$ is open if and only if the set of $\m \in \spec\ R$ such that $R/\m \otimes_R \tilde{M}_A \cong M_A$ is dense in $\spec\ R$ for any one-dimensional domain $R$ and any ${}_R \tilde{M}_A$ corresponding to an $R$-valued point of $\Mod^A_d$ with $R/\m_0 \otimes_R \tilde{M}_A \cong M$ for some $\m_0 \in \spec\ R$.  Thus suppose that we have such an ${}_R \tilde{M}_A$, and assume also that the orbit of $N_B : = M \otimes_A X_B$ is open in $\Mod^B_e$.   Since ${}_A X$ is projective and ${}_R \tilde{M}$ is free, ${}_R \tilde{M} \otimes_A X_B$ is projective over $R$, hence locally free.  As we also have $R/\m_0 \otimes_R (\tilde{M} \otimes_A X_B) \cong N_B$, we can conclude that $R/\m \otimes_R (\tilde{M} \otimes_A X_B) \cong N_B$ for a dense subset of $\m \in \spec\ R$.  Tensoring with $Y$ and using $X \otimes_B Y \cong A \oplus P$, we have $$R/\m \otimes_R \tilde{M}_A \oplus R/\m \otimes_R \tilde{M} \otimes_A P_A \cong M \oplus M \otimes_A P_A$$ for each $\m$ belonging to a dense subset of $\spec\ R$.  Comparing non-projective summands and using the fact that $\dim_k (R/\m \otimes_R \tilde{M}_A) = d = \dim_k M$, we can conclude that $R/\m \otimes_R \tilde{M}_A \cong M_A$ for a dense subset of $\m \in \spec\ R$.  Hence $M$ is rigid and $\Or(M)$ is open.

The converse is proved similarly, as $M \otimes_A X_B \cong N \oplus Q_0$ for some projective $Q_0$ and some $N_B$ with $N \otimes_B Y_A \cong M \oplus P_0$ with $P_0$ projective.  $\Box$\\

We can now prove our main result.

\vspace{5mm}
\noindent
{\it Proof of Theorem 1.}    Since the dimensions of the projective modules appearing in a minimal projective resolution of $M$ are bounded, we see that the dimensions of the syzygies of $M$ assume only finitely many values.  Hence there is some $d$ with $\dim_k \Omega^i M = d$ for infinitely many $i \geq  0$.  Since the module variety $\Mod^A_d$ of $d$-dimensional $A$-modules has a finite number of irreducible components, there is an irreducible component $C$ containing the points corresponding to $\Omega^i M$ for infinitely many $i$.  By the previous result, each $\Omega^i M$, as a summand of $M \otimes_A \Omega^i_{A^e}(A)$, has an open and dense orbit in $C$.  Thus if $\Omega^i M$ and $\Omega^j M$ belong to $C$ for $i > j \geq 0$, then their $\GL_d(k)$-orbits have a nonempty intersection, which yields $\Omega^j M \cong \Omega^i M$.  Since the syzygy functor is an auto-equivalence of $\stmod A$, we now have $M \cong \Omega^{i-j}M$, and $M$ is periodic.  $\Box$\\

\noindent
{\bf Remark.}  In the language of {\it algebraically rigid modules} over a self-injective algebra, Dade has shown the following \cite{Dade}:
\begin{itemize}
\item For each $d >0$, there are only finitely many algebraically rigid modules of dimension $d$, up to isomorphism (Corollary 1.12).
\item If $M$ is algebraically rigid, then so is $\Omega^i M$ for all $i \in \mathbb{Z}$ (Proposition 3.12).
\item Any simple module is algebraically rigid  (Proposition 3.4).
\end{itemize}
The argument above thus shows that any algebraically rigid module of complexity one is periodic, and hence establishes Corollary 2 as well.\\

\section{Algebras of complexity one}

The complexity of a self-injective algebra $A$ is defined as $$\cx(A) = \sup\ \{\cx(M)\ |\ M \in \rmod A\}= \sup\ \{\cx(S)\ |\ S_A \ \mbox{simple}\}.$$
Using \cite{HCFDA} it is easy to see that $\cx(A) = \cx_{A^e}(A)$ where $A^e = A^{\op} \otimes_k A$ is the enveloping algebra of $A$.  Moreover, we say that $A$ is a {\it periodic algebra} if $A$ is a periodic $A^e$-module.  To the best of our knowledge, there are no known examples of self-injective algebras of complexity one which are not periodic.  In this section, we investigate whether the above methods could be applied to show that a self-injective algebra of complexity one is periodic.

From Corollary 2 it follows that over an indecomposable self-injective algebra $A$ of complexity one every simple module is periodic.  By Theorem 1.4 of \cite{GrSnSo}, we know that this latter condition is equivalent to  $\Omega^n_{A^e}(A) \cong {}_1 A_{\sigma}$ for some $\sigma \in \aut_k(A)$  and some $n \geq 1$, where ${}_1 A_{\sigma}$ denotes the {\it twisted} $(A,A)$-bimodule $A$ with the right action twisted by $\sigma$: $a\cdot x \cdot b = axb^{\sigma}$.  In this case, $A$ is periodic if and only if $\sigma$ has finite order in $\out_k(A)$.

By Theorem 1, to show that a self-injective algebra of complexity one is periodic, it would suffice  to show that the bimodule $A$ has an open orbit in $\Mod^{A^e}_d$.  Unfortunately, this condition appears to fail even in the simplest cases.  For consider the (periodic) algebra $A = k[x]/(x^2)$, which has a one-parameter family of outer automorphisms $\sigma_t$ defined by $x^{\sigma_t} = tx$ for $t \in k^*$.  It follows that the $A^e$-modules ${}_1 A_{\sigma_t}$ provide a one-parameter family of non-isomorphic modules in the same irreducible component of $\Mod^{A^e}_2$, and consequently $\Or({}_1 A_1)$ is not open.  The following significantly weaker statement is perhaps the best we can expect to hold in general regarding the geometry of the bimodule ${}_1A_1$.

\begin{propos} For a self-injective algebra $A$, the $A^e$-module $A$ is not a proper degeneration of any other $A^e$-module.\end{propos}

\noindent
{\it Proof.} Suppose that $M \dgn A$ for an $(A,A)$-bimodule $M$.  Then there is a short exact sequence $\ses{U}{U\oplus M}{A}{}{}$ of bimodules by \cite{Zwa}.  Clearly this sequence splits when considered as a sequence of either left or right $A$-modules, yielding $M_A \cong A_A$ and ${}_A M \cong {}_A A$.  We also see that $S \otimes_A M_A \cong S_A$ for all simple modules $S_A$ (and similarly for all simple left modules).  Using arguments from the proof of Theorem 1.4 in \cite{GrSnSo}, we see that $M$ is isomorphic to a twisted bimodule ${}_1 A_{\sigma}$ for some $\sigma \in \aut_k (A)$.  But $\Endo_{A^e}({}_1 A_{\sigma}) \cong Z(A) \cong \Endo_{A^e}(A)$ implies that the orbits of ${}_1 A_{\sigma}$ and $A$ each have dimension $d^2 - \dim_k Z(A)$ by II.3.6 of \cite{Kraft}.  Thus $A \in \overline{\Or({}_1 A_{\sigma})}$ forces $A \in \Or({}_1 A_{\sigma})$. $\Box$ \\


\begin{thebibliography}{99}

\bibitem{Alp} J. Alperin.  \emph{Periodicity in groups.}  Illinois J. Math.  21 (1977), no. 4, 776--783.

\bibitem{TEG} M. Auslander and I. Reiten.  \emph{On a theorem of E. Green on the dual of the transpose.}  Proc. ICRA V, CMS Conf. Proc. 11 (1991), 53-65.


\bibitem{Bon} K. Bongartz. \emph{On degenerations and extensions of finite-dimensional modules.}  Adv. Math. 121 (1996), no. 2, 24--287.

\bibitem{CBS}  W. Crawley-Boevey and J. Schroer.  \emph{Irreducible components of varieties of modules.}   J. Reine Angew. Math.  553  (2002), 201--220. 

\bibitem{Dade} E. Dade.  \emph{Algebraically rigid modules.} Representation theory, II (Proc. Second Internat. Conf., Carleton Univ., Ottawa, Ont., 1979), pp. 195--215, Lecture Notes in Math., 832, Springer, Berlin, 1980. 

\bibitem{Eis} D. Eisenbud.  \emph{Homological algebra on a complete intersection, with an application to group representations.}  Trans. Amer. Math. Soc.  260 (1980), no. 1, 35--64. 

\bibitem{CAVG} D. Eisenbud.  \emph{Commutative algebra.  With a view toward algebraic geometry.} Graduate Texts in Mathematics, 150. Springer-Verlag, New York, 1995.


\bibitem{ErdSko2} K. Erdmann and A. Skowro\'nski.  \emph{Periodic algebras.}  Trends in Representation Theory and Related Topics.  European Math. Soc., Zurich, 2008.

\bibitem{GashPeev} V. Gasharov and I. Peeva.  \emph{Boundedness versus periodicity over commutative local rings.}  Trans. Amer. Math. Soc. 320 (1990), no. 2, 569-580.

\bibitem{GrSnSo} E. L. Green, N. Snashall, and \O. Solberg. \emph{The Hochschild cohomology ring of a selfinjective algebra of finte representation type.}  Proc. Amer. Math. Soc. 131 (2003), no. 11, 3387--3393.

\bibitem{HCFDA} D. Happel.  \emph{Hochschild cohomology of finite-dimensional algebras.}  S\'eminaire d'Alg\`{e}bre Paul Dubreil et Marie-Paul Malliavin, 39\`{e}me Ann\'ee (Paris, 1987/1988),  108--126, Lecture Notes in Math., 1404, Springer, Berlin, 1989.


\bibitem{Kraft} H. Kraft. \emph{Geometric methods in representation theory.} Representations of algebras (Puebla, 1980),  pp. 180–258, Lecture Notes in Math., 944, Springer, Berlin-New York, 1982. 

\bibitem{LiuSchulz} S. Liu and R. Schulz.  \emph{The existence of bounded infinite $DTr$-orbits.}  Proc. Amer. Math. Soc.  122 (1994), no. 4, 1003-1005.

\bibitem{AMRT} J. Rickard.  \emph{Some recent advances in modular representation theory.}  Algebras and Modules I, CMS Conf. Proc. 23 (1998), 157-178.


\bibitem{Zwa} G. Zwara.  \emph{Degenerations of finite-dimensional modules are given by extensions.}  Composito Mathematica 121 (2000), 205-218.

\end{thebibliography}
\end{document}